\documentclass[11pt,reqno]{amsart}
\usepackage{amsfonts}
\usepackage{amsthm}
\usepackage{amsmath}
\usepackage{amssymb}
\usepackage[dvips]{graphics,graphicx}
\usepackage{caption}
\usepackage[T1]{fontenc}
\usepackage{epsfig}
\usepackage{lmodern}
\usepackage{exscale}
\usepackage[latin1]{inputenc}
\usepackage{xcolor}

\usepackage{fancyhdr}
\usepackage{ifthen}
\usepackage{geometry}
\theoremstyle{plain}
\newtheorem{thm}{Theorem}[section]
\newtheorem{pro}{Proposition}[section]
\newtheorem{lem}{Lemma}[section]

\def\eps{\varepsilon}
\def\R{{\mathbb R}}

\begin{document}

\title[No blowup at zero points of the potential]{Excluding blowup at zero points of the potential\\ by means of Liouville-type theorems}

\author[J.-S. Guo]{Jong-Shenq Guo}
\address{Department of Mathematics, Tamkang University,
151, Yingzhuan Road, Tamsui, New Taipei City 25137, Taiwan}
\email{jsguo@\,mail.tku.edu.tw}

\author[Ph. Souplet]{Phlippe Souplet}
\address{Universit\'e Paris 13, Sorbonne Paris Cit\'e,
Laboratoire Analyse, G\'eom\'etrie et Applications, CNRS (UMR 7539),
93430 Villetaneuse, France}
\email{souplet@\,math.univ-paris13.fr}

\thanks{Date: \today.} 

\thanks{This work was partially supported by the Ministry of Science and Technology
of Taiwan (ROC) under the grant 102-2115-M-032-003-MY3.
 The first author would like to thank the support of LAGA - University of Paris 13 / CNRS, where part of this work was done.
The second author is grateful for the support of the National Center for Theoretical Sciences at Taipei for
his visit to Tamkang University, where part of this work was done.}

\thanks{{\em 2000 Mathematics Subject Classification.} Primary: 35K55, {Secondary: 35B44, 35B53}}

\thanks{{\em Key words and phrases:} blowup, potential, Liouville-type theorem}

\begin{abstract}
We prove a local version of a (global) result of Merle and Zaag about ODE behavior of solutions
near blowup points for subcritical nonlinear heat equations.
As an application, for the equation $u_t= \Delta u+V(x)f(u)$,
we rule out the possibility of blowup at zero points of the potential $V$
for monotone in time solutions when $f(u)\sim u^p$ for large $u$, 
 both in the Sobolev subcritical case and in the radial case.
This solves a problem left open in previous work on the subject.
Suitable Liouville-type theorems play a crucial role in the proofs.
\end{abstract}

\maketitle
\setlength{\baselineskip}{18pt}

\section{Introduction}
\setcounter{equation}{0}

In this paper, we consider the following semilinear heat equation with spatially dependent coefficient
in the nonlinearity:
\begin{equation}
\label{e1}
u_t= \Delta u+V(x)f(u), \quad 0<t<T,\ x\in \Omega.
\end{equation}
In the case when $V$ is a positive constant and $f(u)\sim u^p$ with $p>1$,
the blowup behavior of solutions has received considerable attention
in the past decades and a rich variety of phenomena has been discovered
(see, e.g., the monograph \cite{QSbook} and the references therein).
In the case when the potential is nonnegative and nonconstant, it is a natural question whether or not blowup can occur at
zero points of the potential $V$. Although the answer would intuitively seem to be negative at first sight,
it was surprisingly found in \cite{FT00, GLS10, GS11, GLS13} to be positive or negative depending on the situation
 (see Remark~1.2 for details).

The goal of this paper is twofold:
\smallskip

(i) rule out the possibility of blowup at zero points of the potential $V$
for monotone in time solutions of equation \eqref{e1} when $f(u)\sim u^p$ for large $u$.
\smallskip

(ii) prove a local version of a (global) result of Merle-Zaag \cite{MZ98} about ODE behavior of solutions
near blowup points for subcritical nonlinear heat equations. This result, of independent interest,
 and which seems to be new even in the case $V\equiv 1$, 
will be an essential ingredient for (i).
\smallskip

Let us now state our general assumptions:
\begin{eqnarray}
&&\hbox{$\Omega$ is a domain of $\R^n$,}  \label{hyp0}\\
&&\hbox{$f:[0,\infty)\to [0,\infty)$ is a function of class $C^1$,} \label{hyp1} \\ 
&&\hbox{$\displaystyle\lim_{s\to\infty} s^{-p}f(s)=1$ for some $p>1$,} \label{hyp2} \\
&&\hbox{$|f'(s)|\le C(1+s^{p-1})$,\quad $s\ge 0$,} \label{hyp2b} \\
&&\hbox{$V:\overline\Omega\to [0,\infty)$ is a H\"older continuous function.} \label{hyp3}
\end{eqnarray}
Throughout this article, we set
$$p_S=(n+2)/(n-2)_+,\qquad \alpha=1/(p-1),\qquad \kappa=\alpha^\alpha$$
 and we denote the zeroset of $V$ by
$$\mathcal{V}_0=\bigl\{x\in\overline\Omega;\;V(x)=0\bigr\}.$$

Our first main result rules out the possibility of blowup at zero points of the potential~$V$
for monotone in time solutions of \eqref{e1}, under suitable assumptions.
In fact, the case of the homogeneous Dirichlet problem associated with equation \eqref{e1}
with $f(0)=0$ was completely solved in \cite{GS11}.
The more delicate case $f(0)>0$ was left as an open problem.
We here essentially solve it for subcritical $p$, {under a mild geometric assumption.}
Actually, the result here is formulated in a completely local way, without reference to any boundary conditions.
Here, $x_0$ is said to be a blowup point if $\displaystyle\limsup_{t\to T,\, x\to x_0} u(t,x)=\infty$.

\begin{thm}[] \label{th1}
Assume \eqref{hyp0}-\eqref{hyp3}, with $\Omega$ bounded, $f$ of class $C^2$ and convex.
Let $x_0\in\Omega$ be such that $V(x_0)=0$ and
\begin{equation}
\label{HypVzeroset}
\hbox{ the connected component of $\mathcal{V}_0$ containing $x_0$ does not intersect $\partial\Omega$.}
\end{equation}

(i)  Assume that $p<p_S$ and let $u$ be a nonnegative classical solution of \eqref{e1} such that $u_t\ge 0$.
Then $x_0$ is not a blowup point of $u$.

(ii) Assertion (i) remains valid for any $p>1$, if we assume in addition that $u$ and $V$ are radially symmetric
and $\Omega=B_R$.
\end{thm}

 We stress that the assumption $u_t\ge 0$ cannot be removed in general
(compare Theorem~\ref{th1}(ii) with cases (b) and (c) in Remark~1.2 below).
Also, the nonlinearity $u^p$ in Theorem~\ref{th1}
cannot be replaced with a slowly growing one.
Namely, we show in Proposition~\ref{propGlobalBU} below that when $f(u)=u[\log(1+u)]^a$ with $1<a<2$, 
and for suitable potentials whose zeros satisfy \eqref{HypVzeroset},
there exist solutions such that $u_t\ge 0$ and which blow up at every point of the domain.
On the other hand, it is an interesting open problem what happens if assumption \eqref{HypVzeroset}
 is dropped, for instance if there is a line of zeros of $V$
connecting $x_0$ to $\partial\Omega$. The question seems delicate, especially for $n\ge 2$.

\smallskip

Our next main result is a {\it local version} of a global result of Merle-Zaag \cite{MZ98}.
It asserts that the solution of the subcritical nonlinear heat equation
behaves like the corresponding ODE,
in the sense that the diffusion term becomes asymptotically of smaller order than the reaction term
{wherever the solution is large.}
This result is crucial to our proof of Theorem \ref{th1}.

\begin{thm}[] \label{th2}
Let $\omega\subset\subset D\subset\subset \Omega$.
Assume \eqref{hyp0}-\eqref{hyp3}, $p<p_S$ {{and
\begin{equation}
\label{HypVc0}
V(x)\ge c_0 \quad\hbox{in $D$}
\end{equation}
for some $c_0>0$.
Let $u$ be a nonnegative classical solution of \eqref{e1} such that
\begin{equation}
\label{HypTypeI}
u(t,x)\le M(T-t)^{-\alpha} \quad\hbox{in $(T/4,T)\times D$}
\end{equation}
for some $M>0$.}}
Then for {{each $\eps>0$ there}} exists a constant $C_\eps>0$ such that
\begin{equation}
\label{ODEbehav}
\bigl|u_t-V(x)f(u)\bigr|\le \eps u^p+C_\eps  \quad\hbox{in $[T/2,T)\times \omega$.}
\end{equation}
\end{thm}

\smallskip\noindent
{\bf Remark 1.1.} {\it 
(a) Theorem \ref{th2} seems to be new even for $V\equiv 1$.
In the result of \cite{MZ98}, it is assumed that $u$ satisfies the Dirichlet boundary conditions on $\partial\Omega$,
that $\Omega$ is convex, $D=\omega=\Omega$ and $V\equiv 1$.
The case of nonconstant $V$ is mentioned in~\cite[p.~142]{MZ98},
 but (cf.~\cite{Z15}) it is implicitly assumed there that $V$ is positive on $\overline\Omega$.
On the other hand, in \cite{Miz03}, the same global result as in \cite{MZ98} is obtained
for the Neumann problem, without convexity assumption on $\Omega$.

(b) The assumption $p<p_S$ is essentially optimal. Indeed, for $p_S<p<p_L:=1+6/(n-10)_+$ there exist
self-similar, positive classical solutions of $u_t= \Delta u+u^p$, of the form
$u(t,x)=(T-t)^{-\alpha}U(|x|/\sqrt{T-t})$ in $\R^n\times(0,T)$,
with $U$ bounded and $\Delta U(0)<0$ (see \cite{Lep1,BQ, Lep2, Miz09}).
In particular, these solutions satisfy \eqref{HypTypeI} and
$$-\Delta u(t,0)=\theta u^p(t,0)\to \infty,\quad\hbox{ as $t\to T$,}$$
for some $\theta\in (0,1)$,
hence \eqref{ODEbehav} is violated.

(c) Let $D$ be an annulus, $D=\{x\in\R^n:\ r_1<|x|<r_2\}$ with $r_2>r_1>0$,
and assume that $u$ and $V$ are radially symmetric. Then 
Theorem 1.2 remains true for all $p>1$.
Moreover, estimate \eqref{HypTypeI} is also true for all $p>1$ (cf. \cite{PQS07}).
}

\bigskip

The proof of Theorem \ref{th2} is done by rescaling arguments, relying on a Liouville type theorem
of \cite{MZ98} for ancient solutions of $u_t-\Delta u=u^p$.
In this sense it follows the scheme of proof of the corresponding global result in \cite{MZ98}.
However, the proofs in both \cite{MZ98} and \cite{Miz03}
 make extensive use of the so-called weighted energy of Giga and Kohn \cite{GK87},
which requires working with prescribed boundary conditions.
This tool is required in order to guarantee an upper type I estimate, as well as to avoid degeneracy of blowup.
In the local case, we are here able to avoid any energy argument and to
replace this ingredient by a different nondegeneracy property
which is of purely local nature (see Proposition~\ref{prop1}
and Lemma \ref{lem2} below).\footnote{Actually by the proof of Theorem 1.2,
we see that the result remains true if $V$ also depends on $t$ (with $V(t,x)$ H\"older continuous
 satisfying \eqref{HypVc0} in $[T/2,T)\times D$), in which case no energy structure is available in general, even under prescribed boundary conditions.}

As for the local type I estimate (cf. assumption \eqref{HypTypeI} in Theorem \ref{th2}),
it is guaranteed by the following known result, 
which is a consequence of a different Liouville-type Theorem from \cite{PQS07} and \cite{Q16}
(see Theorem 3.1 and Remark 4.3(b) in \cite{PQS07}; note that the proof is given there for constant $V$, but it carries
over with straightforward changes).

\medskip

\noindent{\bf Theorem A.} {\it
Assume \eqref{hyp0}-\eqref{hyp3}, $p<p_S$, $D\subset\subset \Omega$,
\eqref{HypVc0} for some $c_0>0$, and let $u$ be a nonnegative classical solution of \eqref{e1}.
Assume in addition that one of the following conditions holds:
$$\hbox{{\rm (A)} $u_t\ge 0$, \quad {\rm (B)} $n\le 2$, \quad {\rm (C)} $n\ge 3$ and $p<n(n+2)/(n-1)^2$.}$$
Then there exists $M>0$ such that estimate \eqref{HypTypeI} holds.
}

\smallskip

\smallskip\noindent
{\bf Remark 1.2.} {\it 
Let us summarize the previously known results about blowup or non blowup at zero points of the potential 
for positive solutions of the equation
$$
u_t= \Delta u+V(x) u^p, \quad 0<t<T,\ x\in \Omega,
$$
with $p>1$ and homogeneous Dirichlet boundary conditions (if $\Omega\neq\R^n$).

(a) $\Omega=B_R$, $V(x)=|x|^\sigma$, $\sigma>0$, $n\ge 3$, $p<1+2\sigma/(n-1)$ (with equality permitted if $n=3$).
If $u$ is symmetric, then $0$ is not a blowup point~\cite{GLS10};

(b) $\Omega=B_R$, $V(x)=|x|^\sigma$, $\sigma>0$, $n=3$, $p>p_S+2\sigma$.
Then there exist symmetric solutions such that $0$ is a blowup point~\cite{GS11};

(c) $\Omega=\R^n$, $V(x)=|x|^\sigma$, $\sigma>0$, $3\le n\le 10+2\sigma$, $p>p_S+2\sigma/(n-2)$ (or an additional restriction for $n>10+2\sigma$).
Then there exist symmetric (backward self-similar) solutions such that $0$ is a blowup point~\cite{FT00};

(d) $\Omega$ bounded, $0\le V$ H\"older continuous in $\overline\Omega$.
If $u$ is nondecreasing in time,
then blowup cannot occur at any zero point of $V$~\cite{GS11}.

Moreover, in examples (b) and (d), the blowup is of type II.
Finally, for the case of equation \eqref{e1} with $f(0)>0$,
specifically $f(u)=(1+u)^p$ or $f(u)=e^u$, blowup at zero points of $V$ was excluded in \cite{GLS13} under the assumption that
$u$ is nondecreasing in time and the blowup set of $u$ is a compact subset of $\Omega$.
However, for the Dirichlet problem, the latter condition is only known to hold when
$\Omega$ is convex and the potential $V(x)$ is monotonically decreasing near the boundary.}

\medskip

For results on other aspects of equation~\eqref{e1} with variable potential and power-like nonlinearity,
we refer to, e.g., \cite{Pin, QS, LQ, DLY, CER, ABKQ, Xi, Fo}.
For other applications of Merle-Zaag's Liouville theorem and ODE behavior
in the study of the blowup set in the case $V\equiv 1$, see \cite{Z06, FI} and the references therein.
\smallskip

Let us now explain the difficulties in ruling out the possibility of blowup at zero points of the potential $V$,
and the new ideas to overcome these difficulties.
 It is already known \cite{GS11,GLS13} that {\it type~I} blowup cannot occur at a zero point of the potential.
This follows from a local comparison with a suitable self-similar supersolution (recalled in Lemma~\ref{lemomega3} below).
Consequently, we are reduced to proving an a priori type I estimate for the solution $u$ near a blowup point $x_0$.
A classical way to derive such an estimate is to apply the maximum principle to an auxiliary
function of the form $J=u_t-\eps f(u)$, so as to show $J\ge 0$ (cf.~\cite{FM85}).
When $f(0)=0$,  for the Dirichlet problem, one can work on the whole domain $\Omega$,
taking advantage of the fact that $J=0$ on $\partial\Omega$.
When $f(0)>0$ (or without prescribed boundary conditions), this is no longer possible.
When the blowup set is compact (cf. Remark~1.2), this can be fixed by working on
a subdomain of $\Omega$ (cf.~\cite{GLS13}).
For the corresponding quenching problem, with $f(u)=(1-u)^{-p}$, 
the possible lack of compactness of the blowup set, pointed as an obstacle in \cite{GPW, EGG}, was
overcome in \cite{GS15} by introducing a modified functional $J=u_t-\eps a(x)f(u)$,
where $a(x)$ is a suitable function vanishing on the boundary.
However this construction does not seem to work for the blowup problem,
and we here use a completely different idea. Namely by using a Liouville type theorem of Merle and Zaag \cite{MZ98}
(see Theorem~B in section 3 below),
we can show (cf. Theorem \ref{th2}) that, near any point where $V$ is nonzero and $u$ is
large enough, $u_t$ is of same order as $u^p$, hence $J\ge 0$.
If $V(x_0)=0$ and the connected component of $\mathcal{V}_0$ containing $x_0$ does not intersect $\partial\Omega$,
this allows one to apply the maximum principle to $J$ on a suitable subdomain containing $x_0$,
 and on the boundary of which $V$ is positive.

\smallskip

The outline of the paper is as follows. In Section~2, we prove a nondegeneracy result
which is one of the ingredients of the proof of Theorem \ref{th2}. Then Theorem \ref{th2} is proved in Section~3.
Section~4 is devoted to the proof of Theorem~\ref{th1}.
Finally, in Section~5, we consider the case of weak nonlinearities, for which 
we prove blowup at zero points of the potential, and also obtain additional information on the blowup behavior
which stands in contrast with the case of power nonlinearities with a potential.

\bigskip

\section{A nondegeneracy result}
\setcounter{equation}{0}

As an ingredient to the proof of Theorem \ref{th2}, we will prove the following nondegeneracy property,
valid for any $p>1$.
It extends a result of Giga and Kohn (see Theorem 2.1 in \cite{GK89} and cf. also Theorem 25.3 in \cite{QSbook}).
Namely,  condition \eqref{HypIneqParabBound} below involves
 the sharp constant $\kappa$, instead of a small number $\eps>0$ in \cite{GK89}.

\medskip

\begin{pro}[] \label{prop1}
Let $p>1$, $r,T,A>0$, $x_0\in \R^n$ and $u$ be a positive classical solution of the inequality
\begin{equation}
\label{HypIneqParab2c}
u_t-\Delta u\le Au^p\quad\hbox{ in $(0,T)\times B(x_0,r).$}
\end{equation}
Assume that there exist $k\in (0,\kappa)$ and $\delta\in (0,T)$ such that
\begin{equation}
\label{HypIneqParabBound}
(T-t)^\alpha u(t,x) \le kA^{-\alpha}\quad\hbox{ in $[T-\delta,T)\times B(x_0,r)$}.
\end{equation}
Then there exists a constant $L=L(p,n,A,k,r,\delta)>0$ such that
$$u(t,x) \le L\quad\hbox{ in $(T-\delta,T)\times B(x_0,r/4)$}.$$
\end{pro}

\medskip

{\it Proof.}
Assume $x_0=0$ without loss of generality.
For any $\sigma\in (0,2)$ and $R\in (0,r]$, we may find $\phi_R\in C^2(\R^n)$ such that $0\le\phi_R\le 1$,
\begin{equation}
\label{CutOff1}
\phi_R(x)=0\ \hbox{ for $|x|\geq 2R/3$},
\qquad\phi_R(x)=1\ \hbox{ for $|x|\leq R/2$},
\end{equation}
and
\begin{equation}
\label{CutOff2}
|\nabla\phi_R|^2+|\Delta \phi_R^2|\leq C(R,n)\phi_R^\sigma.
\end{equation}
Indeed, this function can be constructed as follows. We fix a nonincreasing function $h\in C^3([1/2,a])$ with $a=2/3$,
 such that $h(1/2)=1$, $h'(1/2)=h''(1/2)=0$,
$h(a)=h'(a)=h''(a)=0$ and $h'''(a)=-6$. We then let
$$
\psi(s)=
\begin{cases}
1,& s\le 1/2 \\
(h(s))^l,& 1/2<s<a\\
0,& s\ge a,
\end{cases}
$$
with $l\ge 2$ an integer. We note that $h''(s)\sim 6(a-s)$,
$h'(s)\sim -3(a-s)^2$ and $h(s)\sim (a-s)^3$, as $s\to a_-$, so that
\begin{equation}
\label{CutOff2a}
\psi(s)\sim (a-s)^{3l},\quad \psi'(s)\sim c_1(a-s)^{3l-1},\quad \psi''(s)\sim c_2(a-s)^{3l-2},
\end{equation}
as $s\to a_-$.
Finally setting $\phi_R(x)=\psi(|x|/R)$, we see that $\phi_R$ has the desired properties
(in particular, \eqref{CutOff2} follows from \eqref{CutOff2a} if $l$ is large enough).

Let $\eps>0$ and put
$$v=v_{\eps,R}=u^{1+\eps}\phi^2, \quad\hbox{with $\phi=\phi_R.$}
$$
For $(t,x)\in (T-\delta,T)\times B_R$, we have
$$\begin{aligned}
v_t-\Delta v
&=(1+\eps)u^\eps u_t\phi^2-(1+\eps)\phi^2(u^\eps\Delta u+\eps u^{\eps-1}|\nabla u|^2)\\
&-4(1+\eps)u^\eps\phi\nabla u\cdot\nabla\phi-u^{1+\eps}\Delta \phi^2.
\end{aligned}$$
Using \eqref{HypIneqParab2c} and $4|u^\eps\phi\nabla u\cdot\nabla\phi|
\leq \eps u^{\eps-1}\phi^2|\nabla u|^2+4\eps^{-1}u^{1+\eps}|\nabla\phi|^2$, we obtain
$$v_t-\Delta v
\leq (1+\eps)Au^{p+\eps}\phi^2+(1+\eps)u^{1+\eps}(4\eps^{-1}|\nabla\phi|^2+|\Delta \phi^2|).$$
Moreover, by \eqref{CutOff2} with $\sigma=2(1+\eps)/(p+\eps)$ and Young's inequality, we get
$$\begin{aligned}
&(1+\eps)u^{1+\eps}(4\eps^{-1}|\nabla\phi|^2+|\Delta \phi^2|) \\
&\qquad\qquad=(1+\eps)u^{1+\eps}\phi^{2(1+\eps)/(p+\eps)}\phi^{-2(1+\eps)/(p+\eps)}(4\eps^{-1}|\nabla\phi|^2+|\Delta \phi^2|)\\
&\qquad\qquad \le \eps Au^{p+\eps}\phi^2+B
\end{aligned}$$
where $B=B(\eps,A,p,n,r)>0$
(the computation is valid at any point such that $\phi>0$, but the conclusion is also true where $\phi=0$, by \eqref{CutOff2}).
It follows that
\begin{equation}
\label{HypIneqParab2b}
v_t-\Delta v \leq (1+2\eps)Au^{p-1}v+B\quad\hbox{in $(T-\delta,T)\times B_R.$}
\end{equation}

Set $m:=(1+2\eps)k^{p-1}$.
Since $k<\kappa=\alpha^\alpha$, we may choose $\eps=\eps(p,k)>0$ sufficiently small, so that $m<\alpha$.
We first use the choice $R=r$.
Using \eqref{HypIneqParab2c} and assumption \eqref{HypIneqParabBound}, we see that {{$v=v_{\eps,r}$}} satisfies
$$
v_t-\Delta v\leq m(T-t)^{-1}v+B\quad\hbox{in $(T-\delta,T)\times B_r.$}
$$
For $K>0$, setting
$$\overline v=\overline v(t):=K(T-t)^{-m}-B (m+1)^{-1}(T-t),\quad T-\delta\leq t<T,$$
we have
$$\overline v_t
=Km(T-t)^{-m-1}+B (m+1)^{-1}=m(T-t)^{-1}\overline v+{{B}}\quad\hbox{in $(T-\delta,T)$}.$$
Moreover, taking
$$K\ge K(\delta,A,p,n,k,r):=B{{(m+1)}} \delta^{m+1}+
(kA^{-\alpha})^{1+\eps}\delta^{m-\alpha(1+\eps)},$$
it follows {{from \eqref{HypIneqParabBound}}} that
$$\overline v(T-\delta)=K\delta^{-m}-B (m+1)^{-1}\delta\ge
\bigl[k(A\delta)^{-\alpha}\bigr]^{1+\eps}\ge v(T-\delta,x),\quad x\in B_r.$$
Since $v=0$ on $(T-\delta,T)\times \partial B_r$, we then deduce from the
comparison principle that $v\leq \overline v$
in $[T-\delta,T)\times  B_r$, hence
\begin{equation}
\label{HypIneqParab3}
u^{1+\eps}\leq K(T-t)^{-m}\quad\hbox{in $(T-\delta,T)\times B_{r/2}$}.
\end{equation}

We next use the choice $R=r/2$, with $\eps$ as above.
Going back to \eqref{HypIneqParab2b} and using \eqref{HypIneqParab3}, we see that $v=v_{\eps,r/2}$ satisfies
$$v_t-\Delta v \leq K_1(T-t)^{-\gamma}v+\tilde B
\quad\hbox{in $(T-\delta,T)\times B_{r/2}$},$$
with
$$\gamma:=\frac{(p-1)m}{1+\eps}<1,\; K_1:=(1+2\eps)AK^{(p-1)/(1+\eps)},\;\tilde B=\tilde B(\eps,A,p,n,r)>0.$$
Setting $K_2=2K_1/(1-\gamma)$ and
$$\overline w=\overline w(t):=K_3\exp[-K_2(T-t)^{1-\gamma}],\quad T-\delta\leq t<T,$$
with $K_3=K_3(\delta,A,p,n,k,r)>0$ sufficiently large, we see that
$$\overline w_t
=2K_1(T-t)^{-\gamma}\overline w
\ge K_1(T-t)^{-\gamma}\overline w+\tilde B
\quad\hbox{in $(T-\delta,T)$},$$
as well as $\overline w(T-\delta)\ge v(T-\delta,x)$ in $B_{r/2}$.
We then deduce from the comparison principle that $v\leq \overline w$
in $[T-\delta,T)\times  B_{r/2}$, hence
$$u^{1+\eps}\leq K_3\quad\hbox{in $(T-\delta,T)\times B_{r/4}$}.$$
This concludes the proof.
\qed

\section{Proof of Theorem~\ref{th2}}
\setcounter{equation}{0}

The proof of Theorem~\ref{th2} relies on the following Liouville type result,
which is due to Merle and Zaag \cite{MZ98}.

\medskip
\noindent{\bf Theorem B.} {\it 
Let $1<p<p_S$, let $T\ge 0$ and let $z$ be a nonnegative classical solution of
$$z_t-\Delta z=z^p\quad\hbox{in $(-\infty,T)\times \R^n$,}$$
such that
$$\sup_{(t,x)\in (-\infty,T)\times \R^n} (T-t)^\alpha z(t,x)<\infty.$$
Then $z$ is independent of $x$.
}

\medskip

For the proof of  Theorem~\ref{th2}, we need the following two lemmas.
They are consequences of Theorem~B and of Proposition~\ref{prop1}, respectively.

\begin{lem}[] \label{lem1}
Assume \eqref{hyp0}-\eqref{hyp3} and
\begin{equation}
\label{HypVLem}
V(x)\ge c_0 \quad\hbox{in $D$}
\end{equation}
for some $c_0>0$ and let $\omega\subset\subset D\ \subset\subset\Omega$.
Let $u$ be a nonnegative classical solution $u$ of~\eqref{e1}.

{{\rm{(i)}}}  Assume that $p<p_S$ and that $u$ satisfies
\begin{equation}
\label{HypTypeILem}
u(t,x)\le M(T-t)^{-\alpha} \quad\hbox{in $(T/4,T)\times D$.}
\end{equation}
Then for {{each $\eps>0$}} there exists a constant $C_\eps>0$ such that
\begin{equation}
\label{ConclLemDelta}
|\Delta u(t,x)|\le \eps (T-t)^{-\alpha-1}+C_\eps  \quad\hbox{in $(T/2,T)\times \omega$}
\end{equation}
and
\begin{equation}
\label{ConclLemGrad}
|\nabla u(t,x)|\le \eps (T-t)^{-\alpha-(1/2)}+C_\eps  \quad\hbox{in $(T/2,T)\times \omega$.}
\end{equation}

{{\rm{(ii)}}  Let $p>1$, let $D$ be an annulus, $D=\{x\in\R^n:\ r_1<|x|<r_2\}$ with $r_2>r_1>0$,
and assume that $u$ and $V$ are radially symmetric.
Then \eqref{ConclLemDelta} and \eqref{ConclLemGrad} are true.}
\end{lem}

\begin{lem}[] \label{lem2}
Let $p>1$, $r,T,A>0$, $x_0\in \R^n$ and $u$ be a positive classical solution of the inequality
\begin{equation}
\label{HypIneqParab4}
u_t-\Delta u\le Au^p\quad\hbox{ in $(0,T)\times B(x_0,r)$.}
\end{equation}
Assume that, for each $\eps>0$, there exists a constant $C_\eps>0$ such that
\begin{equation}
\label{ConclLemDelta2}
\Delta u(t,x) \le \eps (T-t)^{-\alpha-1}+C_\eps  \quad\hbox{in $(T/2,T)\times B(x_0,r).$}
\end{equation}

{{\rm{(i)}}} For any $k\in (0,\kappa)$, there exist constants $\tau_0, K>0$ 
with the following property: If
\begin{equation}
\label{HypTypeILemSmall}
(T-t_0)^\alpha u(t_0,\cdot) \le kA^{-\alpha}\quad\hbox{ in $B(x_0,r)$,\quad for some $t_0\in [(T-\tau_0)_+,T)$}.
\end{equation}
then
$$u(t,x) \le K\quad\hbox{ in $[t_0,T)\times B(x_0,r/4)$}.$$
Here $\tau_0$ depends only on $p,k,A$ and on the constants $C_\eps$,
and $K$ depends only on $p,k,A,n,$ $r,t_0$ and on the constants $C_\eps$.
\smallskip

{{\rm{(ii)}}} If $x_0$ is a blowup point of $u$ (i.e., $\limsup_{t\to T,\, x\to x_0} u(t,x)=\infty$), then
\begin{equation}
\label{HypTypeILemSmall2}
\liminf_{t\to T}(T-t)^\alpha u(t,x_0) \ge \kappa A^{-\alpha}.
\end{equation}
\end{lem}

{\it Proof of Lemma \ref{lem1}.} (i)
Assume the contrary. Then there exist $\eps>0$ and a sequence of points $(t_j,x_j)\in (T/2,T)\times \omega$ such that
\begin{equation}
\label{HypContradDelta0}
|\Delta u(t_j,x_j)|+|\nabla u(t_j,x_j)|^{2(\alpha+1)/(2\alpha+1)}\ge \eps (T-t_j)^{-\alpha-1}+j.
\end{equation}
We may assume that $x_j\to x_\infty\in \overline \omega$.
We also have $t_j\to T$, since otherwise,
$\Delta u(t_j,x_j)$ and $\nabla u(t_j,x_j)$ would be bounded.
We rescale $u$ by setting:
$$v_j(s,y):=\lambda_j^{2\alpha} u(t_j+\lambda_j^2 s,x_j+\lambda_j y),
\quad\hbox{ in $D_j:=\{ -t_j\lambda_j^{-2}<s<1,\ |y|<d\lambda_j^{-1}\}$},$$ 
where $\lambda_j=\sqrt{T-t_j}\to 0$ and $d={\rm dist}(\omega,D^c)>0$.
By a simple computation, we see that
\begin{equation}
\label{Eqnfj}
\partial_s v_j-\Delta_y v_j=V_j(y)f_j(v_j(s,y))
\quad\hbox{ in $D_j$},
\end{equation}
where
$f_j(v):=\lambda_j^{2p/(p-1)} f(\lambda_j^{-2/(p-1)}v)$ and $V_j(y):=V(x_j+\lambda_j y)$.

As a consequence of assumption \eqref{HypTypeILem}, we note that
\begin{equation}
\label{EstimMv}
v_j(s,y)\le M\lambda_j^{2\alpha} (T-t_j-\lambda_j^2 s)^{-\alpha}=M(1-s)^{-\alpha}
\quad\hbox{ in $D_j$}.
\end{equation}
Moreover, we have
$$
\begin{aligned}
|\Delta_y v_j(0,0)|
&+|\nabla_y v_j(0,0)|^{2(\alpha+1)/(2\alpha+1)} \\
&=\lambda_j^{2\alpha+2} \Bigl(|\Delta u(t_j,x_j)|+|\nabla u(t_j,x_j)|^{2(\alpha+1)/(2\alpha+1)}\Bigr) \\
&=(T-t_j)^{\alpha+1} \Bigl(|\Delta u(t_j,x_j)|+|\nabla u(t_j,x_j)|^{2(\alpha+1)/(2\alpha+1)}\Bigr)
\end{aligned}
$$
hence, by \eqref{HypContradDelta0},
\begin{equation}
\label{HypContradDelta2c}
|\Delta_y v_j(0,0)|+|\nabla_y v_j(0,0)|^{2(\alpha+1)/(2\alpha+1)} \ge \eps.
\end{equation}

Note that $f_j(0)=\lambda_j^{2p/(p-1)} f(0)\le 1$ and that, by \eqref{hyp2b},
$$|f'_j(v)|=\lambda_j^{2} |f'(\lambda_j^{-2/(p-1)}v)|\le C\lambda_j^{2} \bigl[1+(\lambda_j^{-2/(p-1)}v)^{p-1}\bigr]
\le C(1+v^{p-1}).$$
As a consequence of \eqref{Eqnfj}, \eqref{EstimMv}, \eqref{hyp3} and interior parabolic estimates,
we then deduce that there exist a function $w\ge 0$ and $\nu\in(0,1)$ such that,
for each compact subset $K$ of $E=(-\infty,1)\times \R^n$,
the sequence $v_j$ converges in $C^{1+(\nu/2),2+\nu}(K)$ to $w$ as $j\to \infty$.
In particular, by \eqref{EstimMv} and \eqref{HypContradDelta2c}, we have
\begin{equation}
\label{EstimMv2}
{{w}}(s,y)\le M(1-s)^{-\alpha} \quad\hbox{in $E$,}
\end{equation}
as well as
\begin{equation}
\label{HypContradDelta3}
|\Delta_y w(0,0)|+|\nabla_y w(0,0)|^{2(\alpha+1)/(2\alpha+1)}\ge \eps.
\end{equation}

Now, for each $(s,y)\in E$, considering separately the cases $w(s,y)>0$ and $w(s,y)=0$
and using assumption \eqref{hyp2}, we see that
$$\lim_j \lambda_j^{2p/(p-1)} f(\lambda_j^{-2/(p-1)}v(s,y))=w^p(s,y).$$
It follows that $w$ is a classical solution of
\begin{equation}
\label{LimitingEq}
\partial_s w-\Delta_y w=V(x_\infty)w^p \quad\hbox{in $E$.}
\end{equation}
Note that $V(x_\infty)>0$ by \eqref{HypVLem}.
By Theorem~B, such a solution with the additional property \eqref{EstimMv2} must
necessarily be spatially homogeneous. This contradicts \eqref{HypContradDelta3}.

\smallskip
(ii) We only sketch the necessary changes.
Since $u$ is radial, setting $\rho=|x|$, $\rho_j=|x_j|$ and $\rho=\rho_j+\lambda_j y$, equation \eqref{Eqnfj} can be written as
$$\partial_s v_j-\partial^2_y v_j=V_j(y)f_j(v_j(s,y))+\frac{n-1}{\rho_j+\lambda_jy}\lambda_j^{2\alpha+2}u_\rho(t_j+\lambda_j^2 s,\rho_j+\lambda_j y).$$
On the other hand, in the radial case in an annulus, by \cite{PQS07}, estimate \eqref{HypTypeILem} is true and moreover
$$|u_\rho(t,\rho)|\le M_1(T-t)^{-\alpha-(1/2)}\quad\hbox{ in $(T/4,T)\times (r_1,r_2)$},$$
for some $M_1>0$. Therefore, since $\lambda_j=\sqrt{T-t_j}$, we have
$$\begin{aligned}
\frac{n-1}{\rho_j+\lambda_jy}\lambda_j^{2\alpha+2}|u_\rho(t_j+\lambda_j^2 s,\rho_j+\lambda_j y)|
&\le \frac{n-1}{r_1}(T-t_j)^{\alpha+1}M_1(T-t_j-\lambda_j^2 s)^{-\alpha-(1/2)} \\
&\le \frac{n-1}{r_1}(T-t_j)^{1/2}M_1(1-s)^{-\alpha-(1/2)}
\end{aligned}
$$
in $D_j$. Consequently, the limiting equation \eqref{LimitingEq} becomes
$$\partial_s w-w_{yy}=V(x_\infty)w^p \quad\hbox{in $(-\infty,1)\times\R$},$$
so that we can conclude as before.
\qed

\medskip

{\it Proof of Lemma \ref{lem2}.}
(i) Let $\eps>0$. By our assumptions, we have
$$u_t\le Au^p+\eps (T-t)^{-\alpha-1}+C_\eps
\quad\hbox{ in $(T/2,T)\times B(x_0,r)$.}$$
Since $k\in (0,\kappa)$, we may choose $B$ so that
\begin{equation}
\label{ChoiceB}
kA^{-\alpha}<B<\kappa A^{-\alpha}.
\end{equation}
Set
$$\phi(t):=B(T-t)^{-\alpha}.$$

We claim that there exist $\eps=\eps(p,k,A)>0$ and $\tau_0=\tau_0(p,k,A,C_\eps)>0$ small, such that
\begin{equation}
\label{IneqPhi1}
\phi'(t)\ge A\phi^p+\eps (T-t)^{-\alpha-1}+C_\eps
\quad\hbox{for all $t\in [(T-\tau_0)_+,T)$.}
\end{equation}
Indeed, since $p\alpha=\alpha+1$, \eqref{IneqPhi1} is equivalent to
$$B\alpha-AB^p-\eps \ge C_\eps (T-t)^{p\alpha}\quad\hbox{for all $t\in [(T-\tau_0)_+,T)$.}$$
Since $B\alpha-AB^p>0$ by \eqref{ChoiceB}, we see that this is true
if we choose $\eps>0$ small and then {{$\tau_0$}} small, with the dependence specified above.

Now \eqref{ChoiceB} and assumption \eqref{HypTypeILemSmall} guarantee that, for each $x\in B(x_0,r)$,
we have $u(t_0,x)\le \phi(t_0)$.
By \eqref{IneqPhi1} and ODE comparison, it follows that
$u(t,x)\le\phi(t)$ for all $t\in [t_0,T)$. Consequently, we have
$$u(t,x)\le B(T-t)^{-\alpha}\quad\hbox{ in $[t_0,T)\times B(x_0,r)$.}$$
Since $B<\kappa A^{-\alpha}$, Proposition~\ref{prop1} then guarantees the desired bound.

\smallskip
(ii) If \eqref{HypTypeILemSmall2} fails, then there exist $k<\kappa$ and $t_0\in ((T-\tau_0)_+,T)$ such that $(T-t_0)^\alpha u(t_0,x_0)<kA^{-\alpha}$,
where $\tau_0$ is given by assertion (i). By continuity there exists $r>0$ such that
$(T-t_0)^\alpha u(t_0,\cdot)<kA^{-\alpha}$ in $B(x_0,r)$.
By assertion (i), it follows that $x_0$ is not a blowup point.
\qed

\medskip

We are now in a position to give the proof of Theorem~\ref{th2},
by combining Lemmas \ref{lem1}-\ref{lem2} and an appropriate rescaling argument.
\medskip

{\it Proof of Theorem~\ref{th2}}.
Assume for contradiction that there exist $c_1>0$ and a sequence $(t_j,x_j)\in [T/2,T)\times \omega$ such that
\begin{equation}
\label{HypContradDelta2}
|\Delta u(t_j,x_j)|\ge c_1u^p(t_j,x_j)+j.
\end{equation}

{\bf Step 1.} {\it Nondegeneracy at points $x_j$.}
First, it follows from \eqref{ConclLemDelta} in Lemma \ref{lem1}(i) and \eqref{HypContradDelta2} that, for all $\eps>0$,
$$c_1u^p(t_j,x_j)+j\le {{|\Delta u(t_j,x_j)|}}\le \eps (T-t_j)^{-\alpha-1}+C_\eps.$$
Therefore, $t_j\to T$ and
\begin{equation}
\label{HypSmallXj}
u(t_j,x_j)\le \eps_j (T-t_j)^{-\alpha}
\end{equation}
with $\eps_j\to 0$ as $j\to\infty$.
Moreover, we may also assume that $x_j\to x_\infty\in \overline \omega$.
Note that $x_\infty$ is in particular a blowup point (i.e., $\limsup_{t\to T,\, x\to x_\infty} u(t,x)=\infty$),
since otherwise by parabolic regularity, $\Delta u(t_j,x_j)$ would be bounded.

{{We claim that there exists a subsequence of $\{(t_j,x_j)\}$ (not relabeled) and a sequence $\hat t_j\to T$ such that
\begin{equation}
\label{ClaimKappa}
\hat t_j \in (0,t_j)\quad\hbox{and}\quad (T-\hat t_j)^\alpha u(\hat t_j,x_j)= \frac{\kappa}{2}V^{-\alpha}(x_\infty).
\end{equation}
To}} prove the claim, 
in view of \eqref{HypSmallXj}, by continuity, it suffices to show
that, for {{each $\eta>0$ and each}} $j_0\ge 1$ there exist $j\ge j_0$ and
$t\in (T-\eta,t_j)$ such that
$$(T-t)^\alpha u(t,x_j)\ge \frac{\kappa}{2}V^{-\alpha}(x_\infty).$$
If this were false, then there would exist $\eta>0$ and $j_0\ge 1$ such that,
for all $j\ge j_0$ and $t\in (T-\eta,t_j)$, $(T-t)^\alpha u(t,x_j)<\frac{\kappa}{2}V^{-\alpha}(x_\infty)$.
For each given $t\in (T-\eta,T)$, since $t_j\to T$, we would have $t_j>t$ for $j$ sufficiently large,
hence $(T-t)^\alpha u(t,x_\infty)\le\frac{\kappa}{2}V^{-\alpha}(x_\infty)$, by  letting $j\to\infty$.
Using the continuity of $V$ and applying Lemma \ref{lem2}(ii) for some $A>V(x_\infty)$
close to $V(x_\infty)$, we would deduce that $x_\infty$ is not a blowup point, which is a contradiction. This proves the claim.

{\bf Step 2.} {\it Rescaling and convergence to a bounded flat profile.}
We rescale similarly as in the proof of Lemma \ref{lem1}, but now taking $\hat t_j$ as rescaling times.
Namely, we set:
$$v_j(s,y):=\lambda_j^{2\alpha} u(\hat t_j+\lambda_j^2 s,x_j+\lambda_j y)
\quad\hbox{ in $\hat D_j:=\{ -\hat t_j\lambda_j^{-2}<s<1,\ |y|<d\lambda_j^{-1}\}$},$$
where $\lambda_j=\sqrt{T-\hat t_j}\to 0$ and $d={\rm dist}(\omega,D^c)>0$.
We have
\begin{equation}
\label{eqnvj}
\partial_s v_j-\Delta_y v_j=
V(x_j+\lambda_j y)\lambda_j^{2p/(p-1)} f(\lambda_j^{-2/(p-1)}v_j(s,y))
\quad\hbox{ in $\hat D_j$}.
\end{equation}

As a consequence of {{assumption \eqref{HypTypeI}}}, we note that
\begin{equation}
\label{EstimMv3}
v_j(s,y)\le M\lambda_j^{2\alpha} (T-\hat t_j-\lambda_j^2 s)^{-\alpha}=M(1-s)^{-\alpha}
\quad\hbox{ in $\hat D_j$.}
\end{equation}
By \eqref{ConclLemDelta}-\eqref{ConclLemGrad}, for all $\eps>0$, we have
\begin{equation}
\label{EstimMv1}
\begin{aligned}
|\nabla_y v_j(s,y)|
&=\lambda_j^{2\alpha+1} |\nabla u(\hat t_j+\lambda_j^2 s,x_j+\lambda_j y)| \\
&\le\lambda_j^{2\alpha+1}\bigl[\eps (T-\hat t_j-\lambda_j^2 s)^{-\alpha-(1/2)}+C_\eps\bigr] \\
&=\eps(1-s)^{-\alpha-(1/2)}+\lambda_j^{2\alpha+1}C_\eps
\quad\hbox{ in $\hat D_j$}
\end{aligned}
\end{equation}
and
\begin{equation}
\label{EstimMv2b}
\begin{aligned}
|\Delta_y v_j(s,y)|
&=\lambda_j^{2\alpha+2} |\Delta u(\hat t_j+\lambda_j^2 s,x_j+\lambda_j y)|  \\
& \le \lambda_j^{2\alpha+2} \bigl[\eps(T-\hat t_j-\lambda_j^2 s)^{-\alpha-1}+C_\eps\bigr] \\
&=\eps(1-s)^{-\alpha-1}+\lambda_j^{2\alpha+2}C_\eps
\quad\hbox{ in $\hat D_j$.}
\end{aligned}
\end{equation}
Also, letting
$$s_j:=\frac{t_j-\hat t_j}{T-\hat t_j}\in (0,1),$$
by \eqref{HypContradDelta2}, we see that
\begin{equation}
\label{HypContradDelta2b}
|\Delta_y v_j(s_j,0)|=\lambda_j^{2\alpha+2} |\Delta u(t_j,x_j)|\ge c_1(T-\hat t_j)^{p\alpha} u^p(t_j,x_j)=c_1v^p_j(s_j,0).
\end{equation}
Moreover, by \eqref{ClaimKappa}, we have
\begin{equation}
\label{vj00}
v_j(0,0)=\lambda_j^{2\alpha}u(\hat t_j,x_j)= \frac{\kappa}{2}V^{-\alpha}(x_\infty).
\end{equation}
As a consequence of \eqref{EstimMv1} and \eqref{vj00}, it follows that
\begin{equation}
\label{vj0cv}
v_j(0,y)\to \frac{\kappa}{2}V^{-\alpha}(x_\infty),
\ \hbox{as $j\to\infty$,\ uniformly for $y$ in compact subsets of $\R^n$.}
\end{equation}

By \eqref{eqnvj}, \eqref{EstimMv3}, \eqref{hyp3} and interior parabolic estimates,
we then deduce that there exist a function $w\ge 0$ and $\nu\in(0,1)$ such that,
for each compact subset $K$ of $[0,1)\times \R^n$,
the sequence $v_j$ converges in $C^{1+(\nu/2),2+\nu}(K)$ to $w$ as $j\to \infty$.
By \eqref{eqnvj} and \eqref{EstimMv2b}, the function $w$ solves
$$\partial_s w(s,y)=V(x_\infty)w^p(s,y),\quad 0\le s<1,\ y\in \R^n$$
with
$$w(0,y)=\frac{\kappa}{2}V^{-\alpha}(x_\infty).$$
Integrating this ODE, we obtain
\begin{equation}
\label{solODE}
w(s,y)=V^{-\alpha}(x_\infty)\kappa(2^{p-1}-s)^{-\alpha},\quad 0\le s<1,\ y\in \R^n.
\end{equation}

{\bf Step 3.} {\it Uniform regularity and flatness of rescaled solution and conclusion.}
We shall now apply Lemma \ref{lem2}(i).
First, since $w(s,y)\le K_0:=V^{-\alpha}(x_\infty)\kappa(2^{p-1}-1)^{-\alpha}$, 
we deduce from Step~2 that
for any $s_0\in (0,1)$, there exists $j_1\ge1$ such that
\begin{equation}
\label{boundvj}
v_j(s_0,y)\le K_0+1,\quad |y|\le 2,\ j\ge j_1.
\end{equation}
By assumption  \eqref{hyp2}, there exists $C>0$ such that $f(u)\le \frac{3}{2}(u+C)^{{p}}$ for all $u\ge 0$. This along with
 \eqref{eqnvj} and the continuity of $V$ implies that
$z_j:=v_j+1$ satisfies
$$\partial_s z_j-\Delta_y z_j=\partial_s v_j-\Delta_y v_j
\le 2V(x_\infty)\lambda_j^{2p/(p-1)} (\lambda_j^{-2/(p-1)}v_j+C)^p\quad\hbox{ in $\hat D_j$},$$
hence
\begin{equation}
\label{ineqzj}
\partial_s z_j-\Delta_y z_j\le 2V(x_\infty) z_j^p,\quad 0\le s<1,\ |y|\le 2,
\end{equation}
for all $j$ sufficiently large.
Let $\tau_0$ be given by Lemma \ref{lem2}(i) with $A=2V(x_\infty)$, $r=2$, $T=1$, $k=\kappa/2$.
Choosing $s_0\in [(1-\tau_0)_+,1)$ such that
$$(1-s_0)^\alpha(K_0+2)<\frac{\kappa}{2}(2V(x_\infty))^{-\alpha},$$
in view of \eqref{boundvj} and \eqref{ineqzj}, it then follows from Lemma \ref{lem2}(i) that
$$v_j(s,y)\le z_j(s,y)\le K,\quad s_0\le s<1,\ |y|\le 1/2,\ j\ge j_1.$$
Going back to equation \eqref{eqnvj}, and using parabolic estimates, we deduce
that $v_j$ actually converges to $w$ in $C^{1+(\nu/2),2+\nu}([1/2,1)\times B_{1/4})$.
In view of \eqref{solODE}, this implies
$$\liminf_{j\to\infty} v_j(s_j,0)\ge \frac{\kappa}{2}V^{-\alpha}(x_\infty)\quad\hbox{and}\quad \lim_{j\to\infty} \Delta v_j(s_j,0)=0.$$
This contradicts \eqref{HypContradDelta2b} and the proof is completed.

Finally, we note that in the radial case in an annulus (cf.~Remark~1.1(c)),
the above proof remains valid for all $p>1$, using assertion (ii) of Lemma \ref{lem1}.
\qed

\section{Proof of Theorem~\ref{th1}}
\setcounter{equation}{0}

The proof is carried out through a series of lemmas.
In what follows we denote $\{V\le k\}=\{x\in \overline\Omega:\, V(x)\le k\}$.
We begin with the following simple topological lemma.

\begin{lem}[] \label{lemomega}
Let $\Omega$ be a bounded domain of $\R^n$
and let $V:\overline\Omega\to [0,\infty)$ be a continuous function.
Let $x_0\in\Omega$ and denote by $A_0$ the connected component of $\mathcal{V}_0$ containing~$x_0$.
If $A_0\cap \partial\Omega=\emptyset$,
then there exists a subdomain $\Omega_0\subset\subset\Omega$ and $\eta>0$ such
that $x_0\in \Omega_0$ and
\begin{equation*}
V\ge \eta>0 \quad\hbox{on $\partial \Omega_0$.}
\end{equation*}
\end{lem}

{\it Proof.}
For any positive integer $m$, we denote by $A_m$ the connected component of
$\{V\le 1/m\}$ containing~$x_0$.

\smallskip
{\bf Step 1.} We claim that we have
\begin{equation}
\label{ClaimOmega1}
\hbox{$A_0\subset A_p\subset A_m$, for all $p\ge m\ge 1$,}
\end{equation}
\begin{equation}
\label{ClaimOmega2}
\hbox{$A_m$ is closed,}
\end{equation}
\begin{equation}
\label{ClaimOmega3}
\hbox{$A_0=\displaystyle\bigcap_{m\ge 1} A_m$.}
\end{equation}

Assertion \eqref{ClaimOmega1} follows directly from $\mathcal{V}_0\subset \{V\le 1/p\}\subset \{V\le 1/m\}$.

To show \eqref{ClaimOmega2}, first note that $\{V\le 1/m\}$ is closed, due to the continuity of $V$,
hence $\overline {A_m}\subset \{V\le 1/m\}$.
Also, since the closure of a connected set is connected,
$\overline{A_m}$ is connected.
Now, since $A_m\subset \overline{A_m}$ and, by definition,
$A_m$ is the largest connected subset of $\{V\le 1/m\}$ containing $x_0$,
we necessarily have $\overline{A_m}=A_m$, hence \eqref{ClaimOmega2}.

Let us prove \eqref{ClaimOmega3}. By \eqref{ClaimOmega1}, we have $A_0\subset
K:=\displaystyle\bigcap_{m\ge 1} A_m$. By \eqref{ClaimOmega2} and the boundedness of $\Omega$, each $A_m$ is compact.
As a decreasing intersection of compact connected sets, $K$ is thus connected (and compact).
Moreover, we have $K\subset\displaystyle\bigcap_{m\ge 1} \{V\le 1/m\}=\mathcal{V}_0$.
Then, since $A_0\subset K$ and, by definition, $A_0$ is the largest connected subset of $\mathcal{V}_0$ containing $x_0$,
we necessarily have $K=A_0$, hence \eqref{ClaimOmega3}.

\smallskip
{\bf Step 2.} Next we claim that there exists $m\ge 1$ such that $A_m\cap \partial\Omega=\emptyset$.
Indeed, otherwise, for each $m\ge 1$, we may find $x_m\in A_m\cap \partial\Omega$.
Since $\partial\Omega$ is compact, up to a subsequence, we may assume $x_m\to y$ for some $y\in\partial\Omega$.
For each $m\ge 1$, by \eqref{ClaimOmega1}, we have $x_p\in A_p\subset A_m$ for all $p\ge m$, hence $y\in A_m$
owing to \eqref{ClaimOmega2}.
Therefore $y\in \displaystyle\bigcap_{m\ge 1} A_m=A_0$ in view of \eqref{ClaimOmega3}.
But this is a contradiction with $A_0\cap \partial\Omega=\emptyset$.

Now let $A'_m$ denote the connected component of $\{V< 1/m\}$ containing $x_0$.
Then $A'_m\subset A_m$ due to $\{V<1/m\}\subset \{V\le 1/m\}$.
Since $A_m\cap \partial\Omega=\emptyset$ and $A_m$ is compact, we have
\begin{equation}
\label{ClaimOmega4}
\hbox{$A'_m\subset\subset \Omega$.}
\end{equation}

We then claim that $A'_m$ is open. Let $x\in A'_m\subset \Omega\cap\{V<1/m\}$.
Since $V$ is continuous, there exists $\rho>0$ such that $B_\rho(x)\subset\{V<1/m\}$.
Therefore, $A'_m\cup B_\rho(x)\subset \{V<1/m\}$ is connected (union of two non-disjoint connected sets)
and contains $x_0$. By the definition of $A'_m$ we deduce that $A'_m\cup B_\rho(x)=A'_m$,
hence $B_\rho(x)\subset A'_m$ and $ A'_m$ is open.

Finally, we observe that $V=1/m>0$ on $\partial A'_m$. Therefore, recalling  \eqref{ClaimOmega4}, we see that
$\Omega_0=A'_m$ has all the required properties. \qed

\medskip
As a consequence of Theorem~\ref{th2} and Lemma~\ref{lemomega}, 
we next prove the following local type I blowup lemma.

\begin{lem}[] \label{lemomega2}
Let the assumptions of Theorem~\ref{th1} be in force.
If $x_0$ is a blowup point, then blowup is of type I near $x_0$. More precisely, there exist $M,\rho>0$ such that
\begin{equation}
\label{Step1TypeI}
u(t,x)\le M(T-t)^{-\alpha}\quad\hbox{in $(0,T)\times B(x_0,\rho)$.}
\end{equation}
\end{lem}

{\it Proof.} Following \cite{FM85}, we set
$$J=u_t-\eps f(u).$$
A standard computation yields, using that $f$ is of class $C^2$ and convex,
\begin{equation}
\label{ComputJ}
\begin{aligned}
J_t-\Delta J
&=(u_t-\Delta u)_t-\eps f'(u)(u_t-\Delta u)-\eps f''(u)|\nabla u|^2 \\
&\leq V(x)f'(u)u_t-\eps V(x)f'(u)f(u)=V(x)f'(u)J\quad\hbox{ in $(0,T)\times \Omega$.}
\end{aligned}
\end{equation}
As a consequence of Lemma~\ref{lemomega},
there exist a subdomain $\Omega_0\subset\subset\Omega$ and $\eta>0$ such
that $x_0\in \Omega_0$ and
\begin{equation}
\label{infV}
V\ge 3\eta>0 \quad\hbox{on $\partial \Omega_0$.}
\end{equation}
Recall that $u_t\ge 0$. Since $u$ blows up, we have $u_t\not\equiv 0$.
Therefore, by the strong maximum principle applied to $u_t$, 
there exists $t_1\in [T/2,T)$ such that
\begin{equation}
\label{estimJ1}
\gamma=\inf_{[t_1,T)\times \overline \Omega_0}u_t>0.
\end{equation}
In particular, assuming $\eps\le \eps_0:=\gamma \Bigl(\sup_{x\in\overline \Omega_0} f(u(t_1,x))\Bigr)^{-1}$, we have
$$J(t_1,\cdot)\ge 0\quad\hbox{ in $\Omega_0$.}$$
As a consequence of Theorem~A, Theorem~\ref{th2} and \eqref{infV},
along with the continuity of $V$, there exists $C>0$ such that
\begin{equation}
\label{estimJ2}
u_t\ge 2\eta f(u)-C \quad\hbox{ in $(0,T)\times \partial \Omega_0$.}
\end{equation}
(In the radial case, cf.~assertion (ii) of Theorem~\ref{th1}, we use Remark~1.1(c).)
Assume $\eps<\min(\eps_0,\eta,\gamma\eta/C)$. Combining \eqref{estimJ1} and \eqref{estimJ2}, we obtain
$$J=u_t-\eps f(u)\ge \max\Bigl[\eta f(u)-C,\ \gamma-\eps f(u)\Bigr]\ge 0$$
 in $(0,T)\times \partial \Omega_0$.
 We can then apply the maximum principle to deduce $J\ge 0$ in $[t_1,T)\times  \Omega_0$.
Using \eqref{hyp2}, by integration, estimate \eqref{Step1TypeI} follows on $[t_1,T)$, hence on $(0,T)$.  \qed

\medskip

Blowup at $x_0$ will be finally ruled out by the following lemma, which shows that
type~I blowup cannot occur at a zero point of the potential.

\begin{lem}[] \label{lemomega3}
Let $p>1$, $C>0$ and $V\ge 0$ be a continuous function on $B(x_0,\rho)$ for some $x_0\in\R^n$ and $\rho>0$.
Let $u\ge 0$ be a classical solution of 
$$u_t\le \Delta u+C\,V(x)(1+u)^p\quad \hbox{in $(0,T)\times B(x_0,\rho)$}$$
and assume that $u$ satisfies the type~I estimate \eqref{Step1TypeI}. 
If $V(x_0)=0$, then $x_0$ is not a blowup point, i.e., there exists $r\in (0,\rho)$
such that $u$ is bounded on $(T/2,T)\times B(x_0,r)$.
\end{lem}

This was proved in \cite{GS11}. For completeness, we reproduce the (supersolution) argument of \cite{GS11} as follows.
\medskip

{\it Proof of Lemma \ref{lemomega3}.}
Let $\rho,M$ be the constants in \eqref{Step1TypeI}. Following \cite{GS11}, we introduce the function
$$w(t,x)=\frac{K}{[q(x)+(T-t)]^\alpha},\quad q(x)=\beta \cos^2\left(\frac{\pi|x-x_0|}{2r}\right),$$
where $r\in(0,\rho/2)$ and the constants $\beta\in (0,1)$ and $K>M$ are to be determined later.
Due to \eqref{Step1TypeI}, we have
$$u(t,x)\le w(t,x)\;\mbox{ for $(t,x)\in(0,T)\times \partial B_0$,}\quad B_0:=\{x\mid |x-x_0|<r\}.$$
Clearly, $u(0,x)\le w(0,x)$ for $x\in\bar{B}_0$, if $K$ is chosen sufficiently large.
We set $f(w)=C(1+w)^p$ and compute
$$
w_t-\Delta w-Vf(w)=\left\{1+\Delta q-(\alpha+1)\frac{|\nabla q|^2}{q+(T-t)}\right\}\frac{\alpha K}{[q+(T-t)]^{\alpha+1}}-Vf(w).
$$
Noting that
$$f(w)[q+(T-t)]^{\alpha+1}\le \ 2CK^p$$
for $K$ sufficiently large, it follows that
$$w_t-\Delta w-Vf(w)\ge 0\quad\mbox{in $(0,T)\times B_0$,}$$
provided
\begin{equation}\label{in11}
1+\Delta q(x)-(\alpha+1)\frac{|\nabla q(x)|^2}{q(x)}-\frac{2C}{\alpha} K^{p-1}V(x)\ge 0\quad\mbox{for all $x\in B_0$.}
\end{equation}
Fixing $K$ and choosing $r$ small enough, since $V(x_0)=0$, we have $(2C/\alpha) K^{p-1}V(x)<1/3$ for all $x\in B_0$.
Since, by direct computation,
$|\Delta q|+q^{-1}|\nabla q|^2\le C(r)\beta$,
by choosing $\beta$ small enough,
 inequality \eqref{in11} follows.
Therefore, the comparison principle gives $u\le w$ in $(0,T)\times B_0$ and the assertion is proved. \qed

\medskip
{\it Proof of Theorem~\ref{th1}.}
By \eqref{hyp2}, we have $f(s)\le C(1+s)^p$ for all $s\ge 0$, with some $C>0$.
The theorem is then a direct consequence of Lemmas~\ref{lemomega2} and \ref{lemomega3}.
\qed

\medskip
\smallskip

\section{Blowup at zero points of the potential for weak nonlinearities}
\setcounter{equation}{0}

Consider the problem
\begin{eqnarray}
u_t &=& u_{xx}+V(x) f(u), \quad t>0,\ x\in (-1,1), \label{globalBU1} \\
u(t,-1)&=&u(t,1)\ =\ 0, \quad t>0, \label{globalBU2} \\
u(0,x)&=&u_0(x), \quad x\in (-1,1),\label{globalBU3} 
\end{eqnarray}
where 
\begin{equation}\label{globalBU4} 
f(u)=u[\log(1+u)]^a.
\end{equation}
The following result, announced after Theorem~\ref{th1},
shows that the nonlinearity $u^p$ cannot be replaced with a slowly growing one in Theorem~\ref{th1}.
 Note that our assumptions allow for instance any even, $C^1$ potential $V$ such that 
 $0\le V\le 1$ and  $V=1$ on $[0,1/3]$. Such potential may vanish at 
isolated points and/or on some subintervals of $(1/3,1)$ and all its zeros will satisfy condition \eqref{HypVzeroset}
 whenever $V(1)>0$.
 \medskip
 
\begin{pro}\label{propGlobalBU}
Assume \eqref{globalBU4} with $1<a<2$. 
Let $V,\phi\in C^1([-1,1])$, with $V, \phi\ge 0$, $V(0)>0$, $\phi(0)>0$ and $\phi(1)=0$.
Assume that $V$ and $\phi$ are even and satisfy:
$$
\hbox{$V', \phi'\le 0$\quad on $[0,1/3]$,}
$$
$$
\hbox{$0\le V(x)\le V(1/3)$, \ $0\le \phi(x)\le \phi(1/3)$\quad on $[1/3,1]$.}
$$
Let $u_0=\lambda\phi$ with $\lambda>0$ and denote by $u\ge 0$ the 
unique, maximal classical solution of problem \eqref{globalBU1}-\eqref{globalBU3} and $T=T(u_0)$ its existence time.

(i) Then $T<\infty$ for all $\lambda>0$ sufficiently large and, whenever $T<\infty$, blowup is global, namely
\begin{equation}\label{ConclGlobalBU}
\lim_{t\to T}u(t,x)=\infty\quad\hbox{ for every $x\in (-1,1)$.}
\end{equation}

(ii) Assume moreover that $\phi\in C^2([-1,1])$ and that, for some $r\in (0,1)$,
 $\phi,V>0$ on $[0,r]$ and $\phi_{xx}\ge 0$ on $[r,1]$.
Then, for all $\lambda>0$ sufficiently large, we have $u_t\ge 0$.
\end{pro}

We note that in the case $V=1$, global blowup for problem \eqref{globalBU1}-\eqref{globalBU4} 
and $1<a<2$ was proved in~\cite{Lac} for smooth bounded domains $\Omega\subset\R^n$
and in \cite{SGKM} for $\Omega=(0,\infty)$. For radially symmetric decreasing solutions in $\Omega=\R^n$,
global blowup as well as further qualitative properties of blow-up solutions were obtained in \cite{GV96}.
For the case of nonconstant, possibly vanishing potential, we will need some specific arguments.
On the other hand, we obtain the following information on the blowup behavior, which shows that
blowup remains of ``type~I'' at each point $x\in (-1,1)$,
 in contrast with the case of blowup at zero points of the potential for power nonlinearities
(cf. cases (b) and (c) in Remark~1.2, and Lemma~\ref{lemomega3}).

\begin{pro}\label{propGlobalBU2} 
Let the assumptions of Proposition~\ref{propGlobalBU}(ii) be in force, with $T<\infty$.
Then there exist constants $C_2\ge C_1>0$
such that
\begin{equation}\label{GlobalBUestim}
C_1(1-|x|)\exp\bigl[{C_1(T-t)^{-1/(a-1)}}\bigr] \le u(t,x)\le u(t,0)\le \exp\bigl[{C_2(T-t)^{-1/(a-1)}}\bigr]
\end{equation}
for all $t\in (T/2,T)$ and $x\in (-1,1)$.
\end{pro}

\smallskip

We turn to the proofs of Propositions \ref{propGlobalBU} and \ref{propGlobalBU2}.
Owing to our assumptions, $u(t,\cdot)$ is even for all $t\in (0,T)$,
but it need not be nonincreasing on $[0,1]$, since $u_0$ and $V$ are not assumed to be so. 
However we shall prove that $u$ has some partial monotonicity properties,
which guarantee the persistence of the maximum at the origin,
 a useful fact for the proof of Propositions \ref{propGlobalBU} and \ref{propGlobalBU2}.
These monotonicity properties, which rely on suitable reflection arguments, remain valid for more general problems
with potential as follows. We note that the assumption $L\ge 1/3$ does not seem easy to relax.

\begin{pro}\label{propGlobalBUmonot}
Let $f:[0,\infty)\to [0,\infty)$ be of class $C^1$.
Let $V, u_0\in C^1([-1,1])$, with $V, \phi\ge 0$, $V(0)>0$, $u_0(0)>0$ and $u_0(1)=0$.
Let $L\in [1/3,1)$ and assume that $V$ and $\phi$ are even and satisfy:
\begin{equation}\label{HypVphi1}
\hbox{$V', u'_0\le 0$ \quad on $[0,L]$,}
\end{equation}
\begin{equation}\label{HypVphi2}
\hbox{$0\le V(x)\le V(L)$, \ $0\le u_0(x)\le u_0(L)$\quad on $[L,1]$.}
\end{equation}
Then for any $\tau>0$ and any classical solution $u$ of problem \eqref{globalBU1}-\eqref{globalBU3} on $(0,\tau]$,
we have
\begin{equation}\label{HypVphi1a}
\hbox{$u_x\le 0$ \quad in $(0,\tau]\times [0,L]$,}
\end{equation}
and
\begin{equation}\label{maxat0}
m(t):=\max_{x\in [-1,1]}u(t,x)=u(t,0),\quad 0<t<\tau.
\end{equation}

\end{pro}

{\it Proof.}
Set $J=(-L,1-2L)$. 
We first claim that
\begin{equation}\label{maxat0claim1}
\hbox{$u(t,x)\ge u(t,x+2L)$ \quad for all $t\in [0,\tau)$ and $x\in J$.}
\end{equation}
Consider the function $w(t,x)=u(t,x)-u(t,x+2L)$.
We have $w(t,-L)=0$, since $u(t,\cdot)$ is even, and also $w(t,1-2L)=u(t,1-2L)\ge 0$.
Moreover,  for all $x\in J$, we have $x\in (-L,L)$, due to $L\ge 1/3$, and $x+2L\in (L,1)$. Therefore,
$V(x)\ge V(L)\ge V(x+2L)$ by \eqref{HypVphi1}-\eqref{HypVphi2}, 
and similarly $w(0,x)=u_0(x)-u_0(x+2L)\ge 0$. In particular, $w$ satisfies
$$w_t-w_{xx}=V(x)f(u(t,x))-V(x+2L)f(u(t,x+2L))\ge V(x)\bigl[f(u(t,x))-f(u(t,x+2L))\bigr]$$
in $(0,\tau)\times J$.
Claim \eqref{maxat0claim1} then follows from the maximum principle.

We next prove \eqref{HypVphi1a}.
By \eqref{maxat0claim1}, we have
$u(t,L+y)\le u(t,y-L)=u(t,L-y)$ for all $y\in (0,1-L)$, hence $u_x(t,L)\le 0$ for all $t\in (0,\tau)$.
Next, by \eqref{HypVphi1}, $z=u_x$ satisfies $z(0,x)\le 0$ in $[0,L]$ and
 is a (strong) solution of
$$z_t-z_{xx}=V'(x)f(u)+V(x)f'(u)z\le V(x)f'(u)z\quad\hbox{ in $(0,\tau)\times(0,L)$}.$$
Since also $z(t,0)=0$, property \eqref{HypVphi1a} follows from the maximum principle.

Now, by \eqref{HypVphi1a}, for all $x\in (-L,L)$, we have $u(t,x)\le u(t,0)$.
Moreover, for all $x\in (L,1)$, we have $x-2L\in J\subset (-L,L)$ due to $L\ge 1/3$.
Consequently, by \eqref{maxat0claim1}, we have $u(t,x)\le u(t,x-2L)\le u(t,0)$.
Property \eqref{maxat0}  follows. \qed

\medskip

{\it Proof of Proposition \ref{propGlobalBU}.} (i)

{\bf Step 1.} {\it Finite time blowup and lower blowup estimate.}
The fact that $T<\infty$ for all $\lambda>0$ sufficiently large is a consequence of a standard Kaplan-type argument
(see e.g. \cite[Chapter 17]{QSbook}).
More precisely, taking $\ell\in (0,1)$ and $c>0$ such that $V,\phi\ge c$ on $[-\ell,\ell]$,
one derives a differential inequality for the functional
$\varphi(t)=\int_{-\ell}^\ell u(t,x)\cos(\pi x/2\ell)\, dx$ by using Jensen's inequality ($f$ being convex),
and one concludes $T<\infty$ by using the fact that $\int_2^\infty \frac{ds}{f(s)}<\infty$.

Next, by Proposition \ref{propGlobalBUmonot}, we have $m(t):=\max_{x\in [-1,1]}u(t,x)=u(t,0)$.
Since $u_{xx}(t,0)\le 0$, it follows that
$$m'(t)\le V(0)f(m(t))\le Cm(t)[\log m(t)]^a.$$ 
By integration, using $\lim_{t\to T}m(t)=\infty$, we deduce
\begin{equation}\label{lowerexp}
u(t,0)\ge  \exp\bigl[{C(T-t)^{-1/(a-1)}}\bigr],\quad 0<t<T.
\end{equation}

\smallskip

{\bf Step 2.} {\it Comparison with a linear problem with fast boundary blowup source.}

By the maximum principle, we have $u\ge v$, where $v$ is the solution of the linear problem:
$$\begin{aligned}
v_t &= v_{xx}, \quad 0<t<T,\ x\in (0,1), \\
v(t,0)&=u(t,0), \quad 0<t<T, \\
v(t,1)&=0, \quad 0<t<T, \\
v(0,x)&=u_0(x), \quad x\in (0,1).
\end{aligned}$$
The function $v$ admits the following representation:
$$v(t,x)=\int_0^1 G(t,x;y)u_0(y)\,dy+\int_0^t \frac{\partial G}{\partial y}(t-s,x;0)u(s,0)\,ds,\quad 0<t<T,\ x\in (0,1),$$
where $G$ is the Dirichlet heat kernel of the interval $(0,1)$.
 This problem is studied in detail in \cite{SGKM} when $(0,1)$ is replaced by the half-line $(0,\infty)$, 
taking advantage of the explicit Gaussian heat kernel.
Although this is not available in our case, it is known~\cite{Zh} that $G$ satisfies the following sharp lower estimate:
$$G(t,x;y)\ge c_1\min\Bigl(\frac{\rho(x)\rho(y)}{t},1\Bigr)t^{-1/2}
\exp\bigl[-c_2|x-y|^2/t\bigr],\quad t\in (0,T],\ x,y\in [0,1],$$
for some constants $c_1,c_2>0$ (depending on $T$), where $\rho(x)=\min(x,1-x)$ is the distance to the boundary
 of $(-1,1)$. Consequently, for each $x\in (0,1)$, we have
$$ \frac{\partial G}{\partial y}(t,x;0)=\lim_{y\to 0}\frac{G(t,x;y)}{y}\ge c_1 \rho(x)t^{-3/2}
\exp\bigl[-c_2|x|^2/t\bigr],\quad t\in (0,T].$$
Set $b=1/(a-1)>1$. Using \eqref{lowerexp}, it follows that, for all $t\in (T/2,T)$ and $x\in (0,1)$,
$$\begin{aligned}
v(t,x)
&\ge c_1 \rho(x) \int_0^t (t-s)^{-3/2}\exp\bigl[-c_2|x|^2(t-s)^{-1}\bigr]
\exp\bigl[C(T-s)^{-b}\bigr]\,ds \\
&\ge c_1 \rho(x) \int_{t-(T-t)}^{t-(T-t)/2} (t-s)^{-3/2}\exp\bigl[-c_2|x|^2(t-s)^{-1}\bigr]
\exp\bigl[C(T-s)^{-b}\bigr]\,ds \\
&\ge \frac{c_1}{2} \rho(x)(T-t)^{-1/2}\exp\bigl[-2c_2|x|^2(T-t)^{-1}+2^{-b}C(T-t)^{-b}\bigr].
\end{aligned}$$
 Therefore, since $b>1$, we obtain
\begin{equation}\label{maxat0claim3}
u(t,x)\ge v(t,x)\ge C_1 \rho(x)\exp\bigl[{C_1(T-t)^{-1/(a-1)}}\bigr],
\quad t\in (T/2,T),\ x\in (0,1).
\end{equation}
This, along with \eqref{lowerexp}, guarantees \eqref{ConclGlobalBU}.
\smallskip

(ii) This follows from a standard maximum principle argument (see, e.g., \cite[Section~52.6]{QSbook})
observing that $u_{0,xx}+V(x)u_0[\log(1+u_0)]^a=\lambda\bigl(\phi_{xx}+V(x)\phi[\log(1+\lambda\phi)]^a\bigr)\ge 0$
in $[-1,1]$ for $\lambda>0$ large, in view of our assumptions on $V$ and $\phi$. \qed

\medskip

{\it Proof of Proposition \ref{propGlobalBU2}.}
Set $J=u_t-\eps f(u)$ and note that $J=0$ for $x=\pm 1$.
Since $u_t\ge 0$ by Proposition~\ref{propGlobalBU}(ii) and $u_t\not\equiv 0$ due to $T<\infty$,
 it follows from the Hopf Lemma that $J(t_1,x)\ge 0$ in $(-1,1)$ for some  $t_1\in [T/2,T)$,
if we choose $\eps>0$ sufficiently small.
By the computation in \eqref{ComputJ} (valid for any convex $f\in C^2$) and the maximum principle,
we deduce that $J\ge 0$ in $(t_1,T)\times (-1,1)$. Since $V(0)>0$, we deduce the upper estimate of $u(t,0)$
upon integration,
 whereas the inequality $u(t,x)\le u(t,0)$ is a consequence of~\eqref{maxat0}.

As for the lower estimate of $u(t,x)$ in \eqref{GlobalBUestim}, for $1/3\le |x|<1$ it follows from \eqref{maxat0claim3}.
For $0\le |x|<1/3$, it follows from $u(t,x)\ge u(t,1/3)$, in view of \eqref{HypVphi1a} with $L=1/3$.
\qed

\bibliographystyle{plain}

\begin{thebibliography}{99}


\bibitem{ABKQ}
N. Ackermann, T. Bartsch, P. Kaplick\' y, P. Quittner,
A priori bounds, nodal equilibria and connecting orbits in indefinite superlinear parabolic problems,
{\it Trans. Amer. Math. Soc.} 360 (2008), 3493--3539.

\vskip 2pt
\bibitem{BQ}
C. Budd, Y.-W. Qi,
The existence of bounded solutions of a semilinear elliptic equation,
{\it J. Differential Equations} 82 (1989), 207--218.

\vskip 2pt
\bibitem{CER}
C. Cortazar, M. Elgueta, J.D. Rossi, 
The blow-up problem for a semilinear parabolic equation with a potential,
{\it J. Math. Anal. Appl.} 335 (2007), 418--427.

\vskip 2pt
\bibitem{DLY}
Y. Deng, Y. Li, F. Yang,
On the stability of the positive steady states for a nonhomogeneous semilinear Cauchy problem,
{\it J. Differential Equations} 228 (2006), 507--529.


\vskip 2pt
\bibitem{EGG}
P. Esposito, N. Ghoussoub, Y. Guo, 
Mathematical Analysis of Partial Differential Equations Modeling Electrostatic MEMS, 
Courant Lect. Notes Math. 20, AMS, Providence, RI, 2010.
\vskip 2pt

\bibitem{FT00}
S. Filippas, A. Tertikas,
On similarity solutions of a heat equation with a nonhomogeneous nonlinearity,
{\it J. Differential Equations} 165 (2000), 468--492.

\vskip 2pt
\bibitem{Fo}
J. F\" oldes, 
Liouville theorems, a priori estimates, and blow-up rates for solutions of indefinite superlinear parabolic problems,
{\it Czechoslovak Math. J.} 61 (2011), 169--198.


\vskip 2pt
\bibitem{FM85} A. Friedman, B. McLeod,
Blow-up of Positive Solution of semilinear Heat Equations,
{\it Indiana Univ. Math. J.} 34 (1985), 425--447.


\vskip 2pt
\bibitem{FI} Y. Fujishima, K. Ishige,
Blow-up set for type I blowing up solutions for a semilinear heat equation,
{\it Ann. Inst. H. Poincaré Anal. Non Lin\'eaire} 31 (2014), 231--247. 


\vskip 2pt
\bibitem{GV96} V.A. Galaktionov and J.L. V\'azquez,
Blowup for quasilinear heat equations described by means of nonlinear Hamilton-Jacobi equations,
{\it J. Differential Equations} 127 (1996), 1-40.

\vskip 2pt
\bibitem{GK85} {Y. Giga, R.V. Kohn},
Asymptotically self-similar blow-up of semilinear heat equations,
{\it Comm. Pure Appl. Math.} 38 (1985), 297--319.

\vskip 2pt
\bibitem{GK87} Y. Giga, R.V. Kohn,
Characterizing blowup using similarity variables, 
{\it Indiana Univ. Math. J.} 36 (1987), 1--40.

\vskip 2pt
\bibitem{GK89} Y. Giga, R.V. Kohn,
Nondegeneracy of blowup for semilinear heat equations,
{\it Comm. Pure Appl. Math.} 42 (1989), 845--884.

\vskip 2pt
\bibitem{GLS10} J.-S. Guo, C.-S. Lin, M. Shimojo,
Blow-up behavior for a parabolic equation with spatially dependent coefficient,
{\it Dynam. Systems Appl.} 19 (2010),
415--434.

\vskip 2pt
\bibitem{GLS13} J.-S. Guo, C.-S. Lin, M. Shimojo,
Blow-up for a reaction-diffusion equation with variable coefficient,
{\it Applied Mathematics Letters} 26 (2013), 150--153.

\vskip 2pt
\bibitem{GS11} J.-S. Guo, M. Shimojo,
Blowing up at zero points of potential for an initial boundary value problem,
{\it Comm. Pure Appl. Anal.} 10 (2011), 161--177.

\vskip 2pt
\bibitem{GS15} J.-S. Guo, Ph. Souplet,
No touchdown at zero points of the permittivity profile for the MEMS problem,
{\it SIAM J. Math. Anal.} 47 (2015), 614--625.

\vskip 2pt

\bibitem{GPW}Y. Guo, Z. Pan, M.J. Ward, 
Touchdown and pull-in voltage behavior of a MEMS device with varying dielectric properties,
{\it SIAM J. Appl. Math.} 66 (2005), 309--338.

\vskip 2pt

\bibitem{Lac}
A.A. Lacey, 
Global blow-up of a nonlinear heat equation,
{\it Proc. Roy. Soc. Edinburgh Sect.} A 104 (1986), 161-167.


\vskip 2pt
\bibitem{Lep1} L.A. Lepin,
Countable spectrum of eigenfunctions of a nonlinear heat conduction equation with distributed parameters,
{\it Differentsial'nye Uravneniya} 24 (1988), 1226--1234
(English translation: {\it Differential Equations} 24 (1988), 799--805).

\vskip 2pt
\bibitem{Lep2} L.A. Lepin,
Self-similar solutions of a semilinear heat equation,
{\it Mat. Model.} 2 (1990), 63--74 (in Russian).

\vskip 2pt
\bibitem{LQ}
J. L\'opez-G\'omez, P. Quittner,
Complete and energy blow-up in indefinite superlinear parabolic problems,
{\it Discrete Contin. Dyn. Syst.} 14 (2006), 169--186.


\vskip 2pt
\bibitem{MZ98} F. Merle, H. Zaag,
Optimal estimates for blowup rate and behavior for nonlinear heat equations,
{\it Comm. Pure Appl. Math.} 51 (1998), 139--196.

\vskip 2pt
\bibitem{Miz03} N. Mizoguchi,
Blowup rate of solutions for a semilinear heat equation with the Neumann boundary condition,
{\it J. Differential Equations} 193 (2003), 212--238.

\vskip 2pt
\bibitem{Miz09} N. Mizoguchi,
Nonexistence of backward self-similar blowup solutions to a supercritical semilinear heat equation,
{\it J. Funct. Anal.} 257 (2009), 2911--2937.


\vskip 2pt
\bibitem{Pin}
R.G. Pinsky, 
Existence and nonexistence of global solutions for $u_t = \Delta u + a(x)u^p$ in $R^d$, 
{\it J. Differential Equations} 133 (1997), 152--177.


\vskip 2pt
\bibitem{PQS07} P. Pol\'a\v cik, P. Quittner, Ph. Souplet,
Singularity and decay estimates in superlinear problems via Liouville-type theorems.
Part II: Parabolic equations, 
{\it Indiana Univ. Math. J.} 56 (2007), 879--908.



\vskip 2pt
\bibitem{Q16} P. Quittner, 
Liouville theorems for scaling invariant superlinear parabolic problems with gradient structure,
{\it Math. Ann.} 364 (2016), 269--292.

\vskip 2pt
\bibitem{QS}
P. Quittner, F. Simondon,
A priori bounds and complete blow-up of positive solutions of indefinite superlinear parabolic problems,
{\it J. Math. Anal. Appl.} 304 (2005), 614--631. 


\vskip 2pt
\bibitem{QSbook} P. Quittner, Ph. Souplet,
Superlinear Parabolic Problems Blow-up, Global Existence and Steady States,
Birkh\"auser Verlag AG, Basel Boston Berlin, 2007.

\vskip 2pt

\bibitem{SGKM}
A.A. Samarskii, V.A. Galaktionov, S.P. Kurdyumov, A.P. Mikhailov,
Blow-up in quasilinear parabolic equations
Nauka, Moscow, 1987. English translation: Walter de Gruyter, Berlin, 1995.

\vskip 2pt
\bibitem{Xi}
R. Xing, 
The blow-up rate for positive solutions of indefinite parabolic problems and related Liouville type theorems,
{\it Acta Math. Sin. (Engl. Ser.)} 25 (2009), 503--518. 

\vskip 2pt
\bibitem{Z06} H. Zaag,
Determination of the curvature of the blow-up set and refined singular behavior for a semilinear heat equation,
{\it Duke Math. J.} 133 (2006), 499--525.

\vskip 2pt
\bibitem{Z15} H. Zaag,
Personal communication, 2015.

\vskip 2pt

\bibitem{Zh}
Q.S. Zhang, 
The boundary behavior of heat kernels of Dirichlet Laplacians,
{\it J. Differential Equations} 182 (2002), 416-430.

\end{thebibliography}

\end{document}